\documentclass[11pt,a4paper,fleqn,leqno]{article}
\usepackage[text={388.37646pt,591.5302pt},vmarginratio=4:5]{geometry}

\usepackage[latin1]{inputenc} %% but in fact source is in 7bit
\def\ladate{v1, 13 June 2011; v2, 4 August 2011}

\usepackage[T1]{fontenc}
\usepackage[english]{babel}

\usepackage{amsmath,amsthm,mathtools}
\allowdisplaybreaks
\def\ipi #1,#2\stop{\mathopen(#1\mathord,#2\mathclose)}
\newcommand\ip[1]{\ipi #1\stop}

\newtheorem{theo}{Theorem}

\newtheorem*{theo*}{Theorem}
\newtheorem*{cor*}{Corollary}
\newtheorem*{prop*}{Proposition}
\theoremstyle{definition}
\newtheorem{note}{Note}

\usepackage[pdfencoding=pdfdoc]{hyperref}

\hypersetup{%
colorlinks,%
linkcolor={blue},%
citecolor={red},%
bookmarksopen,%
pdfdisplaydoctitle,%
pdfpagelayout=OneColumn,%
pdfstartview=FitH,%
pdftitle={On the system of the functions zeta(s)/(s-rho)^k},%
pdfauthor={Jean-Fran\c{c}ois Burnol},%
pdfkeywords={Zeta function, Hilbert spaces, Frames},%
pdfsubject={Number theory, Operator Theory}}

\usepackage{bookmark}
\newcommand\jfautoref[1]{\autoref{#1}}
\newcommand\jfautonote[1]{\autoref{#1}}
\usepackage{txfonts}

\renewcommand{\familydefault}{\ttdefault}
\renewcommand{\ttdefault}{txvtt}
\makeatletter

  \begingroup
  \nfss@catcodes
\DeclareFontFamily{T1}{txvtt}{\hyphenchar\font=127%
\fontdimen2\font=0.33333\fontdimen6\font%
\fontdimen3\font=0.16666\fontdimen6\font%
\fontdimen4\font=0.11111\fontdimen6\font}%

\DeclareFontShape{T1}{txvtt}{m}{n}{	%medium
     <->t1xtt%
}{}%
\DeclareFontShape{T1}{txvtt}{m}{sc}{	%cap & small cap
     <->t1xttsc%
}{}%
\DeclareFontShape{T1}{txvtt}{m}{sl}{	%slanted
     <->t1xttsl%
}{}%
\DeclareFontShape{T1}{txvtt}{m}{it}{	%italic
     <->t1xttsl%
}{}%
\DeclareFontShape{T1}{txvtt}{m}{ui}{   	%unslanted italic
     <->t1xttsl%
}{}%
\DeclareFontShape{T1}{txvtt}{bx}{n}{	%bold extended
     <->t1xbtt%
}{}%
\DeclareFontShape{T1}{txvtt}{bx}{sc}{	%bold extended cap & small cap
     <->t1xbttsc%
}{}%
\DeclareFontShape{T1}{txvtt}{bx}{sl}{	%bold extended slanted
     <->t1xbttsl%
}{}%
\DeclareFontShape{T1}{txvtt}{bx}{it}{	%bold extended italic
     <->t1xbttsl%
}{}%
\DeclareFontShape{T1}{txvtt}{bx}{ui}{  	%bold extended unslanted italic
     <->t1xbttsl%
}{}%
\DeclareFontShape{T1}{txvtt}{b}{n}{	%bold
     <->t1xbtt%
}{}%
\DeclareFontShape{T1}{txvtt}{b}{sc}{	%bold cap & small cap
     <->t1xbttsc%
}{}%
\DeclareFontShape{T1}{txvtt}{b}{sl}{	%bold slanted
     <->t1xbttsl%
}{}%
\DeclareFontShape{T1}{txvtt}{b}{it}{   	%bold italic
     <->t1xbttsl%
}{}%
\DeclareFontShape{T1}{txvtt}{b}{ui}{   	%bold unslanted italic
     <->t1xbttsl%
}{}%
  \endgroup

\re@DeclareMathSymbol{\Gamma}{\mathalpha}{lettersA}{0}
\re@DeclareMathSymbol{\Sigma}{\mathalpha}{lettersA}{6}
\let\alpha\alphaup
\let\beta\betaup
\let\gamma\gammaup
\let\delta\deltaup
\let\epsilon\epsilonup
\let\zeta\zetaup
\let\eta\etaup
\let\theta\thetaup
\let\iota\iotaup
\let\kappa\kappaup
\let\lambda\lambdaup
\let\mu\muup
\let\nu\nuup
\let\xi\xiup
\let\pi\piup
\let\rho\rhoup
\let\sigma\sigmaup
\let\tau\tauup
\let\upsilon\upsilonup
\let\phi\phiup
\let\chi\chiup
\let\psi\psiup
\let\omega\omegaup

  \def\m@t@enc{\encodingdefault}
  \def\m@t@fam{\familydefault}
  \def\m@t@ser{\seriesdefault}
  \def\m@t@opsh{\shapedefault}  
  \def\m@t@bold{\bfdefault}
\DeclareSymbolFont{mtoperatorfont}
    {\m@t@enc}{\m@t@fam}{\m@t@ser}{\m@t@opsh}
\DeclareSymbolFontAlphabet{\mathrm}{mtoperatorfont}
\DeclareMathAlphabet{\mathbf}{\m@t@enc}{\m@t@fam}{\m@t@bold}{\m@t@opsh}
  \SetSymbolFont{mtoperatorfont}{bold} {\m@t@enc}
                                       {\m@t@fam}
                                       {\m@t@ser}
                                       {\m@t@opsh}
  \SetMathAlphabet{\mathbf}{bold}{\m@t@enc}
                                 {\m@t@fam}
                                 {\m@t@ser}
                                 {\m@t@opsh}
\def\operator@font{\mathgroup\symmtoperatorfont}

\DeclareMathSymbol{a}{\mathalpha}{mtoperatorfont}{`a}
\DeclareMathSymbol{b}{\mathalpha}{mtoperatorfont}{`b}
\DeclareMathSymbol{c}{\mathalpha}{mtoperatorfont}{`c}
\DeclareMathSymbol{d}{\mathalpha}{mtoperatorfont}{`d}
\DeclareMathSymbol{e}{\mathalpha}{mtoperatorfont}{`e}
\DeclareMathSymbol{f}{\mathalpha}{mtoperatorfont}{`f}
\DeclareMathSymbol{g}{\mathalpha}{mtoperatorfont}{`g}
\DeclareMathSymbol{h}{\mathalpha}{mtoperatorfont}{`h}
\DeclareMathSymbol{i}{\mathalpha}{mtoperatorfont}{`i}
\DeclareMathSymbol{j}{\mathalpha}{mtoperatorfont}{`j}
\DeclareMathSymbol{k}{\mathalpha}{mtoperatorfont}{`k}
\DeclareMathSymbol{l}{\mathalpha}{mtoperatorfont}{`l}
\DeclareMathSymbol{m}{\mathalpha}{mtoperatorfont}{`m}
\DeclareMathSymbol{n}{\mathalpha}{mtoperatorfont}{`n}
\DeclareMathSymbol{o}{\mathalpha}{mtoperatorfont}{`o}
\DeclareMathSymbol{p}{\mathalpha}{mtoperatorfont}{`p}
\DeclareMathSymbol{q}{\mathalpha}{mtoperatorfont}{`q}
\DeclareMathSymbol{r}{\mathalpha}{mtoperatorfont}{`r}
\DeclareMathSymbol{s}{\mathalpha}{mtoperatorfont}{`s}
\DeclareMathSymbol{t}{\mathalpha}{mtoperatorfont}{`t}
\DeclareMathSymbol{u}{\mathalpha}{mtoperatorfont}{`u}
\DeclareMathSymbol{v}{\mathalpha}{mtoperatorfont}{`v}
\DeclareMathSymbol{w}{\mathalpha}{mtoperatorfont}{`w}
\DeclareMathSymbol{x}{\mathalpha}{mtoperatorfont}{`x}
\DeclareMathSymbol{y}{\mathalpha}{mtoperatorfont}{`y}
\DeclareMathSymbol{z}{\mathalpha}{mtoperatorfont}{`z}
\DeclareMathSymbol{A}{\mathalpha}{mtoperatorfont}{`A}
\DeclareMathSymbol{B}{\mathalpha}{mtoperatorfont}{`B}
\DeclareMathSymbol{C}{\mathalpha}{mtoperatorfont}{`C}
\DeclareMathSymbol{D}{\mathalpha}{mtoperatorfont}{`D}
\DeclareMathSymbol{E}{\mathalpha}{mtoperatorfont}{`E}
\DeclareMathSymbol{F}{\mathalpha}{mtoperatorfont}{`F}
\DeclareMathSymbol{G}{\mathalpha}{mtoperatorfont}{`G}
\DeclareMathSymbol{H}{\mathalpha}{mtoperatorfont}{`H}
\DeclareMathSymbol{I}{\mathalpha}{mtoperatorfont}{`I}
\DeclareMathSymbol{J}{\mathalpha}{mtoperatorfont}{`J}
\DeclareMathSymbol{K}{\mathalpha}{mtoperatorfont}{`K}
\DeclareMathSymbol{L}{\mathalpha}{mtoperatorfont}{`L}
\DeclareMathSymbol{M}{\mathalpha}{mtoperatorfont}{`M}
\DeclareMathSymbol{N}{\mathalpha}{mtoperatorfont}{`N}
\DeclareMathSymbol{O}{\mathalpha}{mtoperatorfont}{`O}
\DeclareMathSymbol{P}{\mathalpha}{mtoperatorfont}{`P}
\DeclareMathSymbol{Q}{\mathalpha}{mtoperatorfont}{`Q}
\DeclareMathSymbol{R}{\mathalpha}{mtoperatorfont}{`R}
\DeclareMathSymbol{S}{\mathalpha}{mtoperatorfont}{`S}
\DeclareMathSymbol{T}{\mathalpha}{mtoperatorfont}{`T}
\DeclareMathSymbol{U}{\mathalpha}{mtoperatorfont}{`U}
\DeclareMathSymbol{V}{\mathalpha}{mtoperatorfont}{`V}
\DeclareMathSymbol{W}{\mathalpha}{mtoperatorfont}{`W}
\DeclareMathSymbol{X}{\mathalpha}{mtoperatorfont}{`X}
\DeclareMathSymbol{Y}{\mathalpha}{mtoperatorfont}{`Y}
\DeclareMathSymbol{Z}{\mathalpha}{mtoperatorfont}{`Z}

\DeclareMathSymbol{0}{\mathalpha}{mtoperatorfont}{`0}
\DeclareMathSymbol{1}{\mathalpha}{mtoperatorfont}{`1}
\DeclareMathSymbol{2}{\mathalpha}{mtoperatorfont}{`2}
\DeclareMathSymbol{3}{\mathalpha}{mtoperatorfont}{`3}
\DeclareMathSymbol{4}{\mathalpha}{mtoperatorfont}{`4}
\DeclareMathSymbol{5}{\mathalpha}{mtoperatorfont}{`5}
\DeclareMathSymbol{6}{\mathalpha}{mtoperatorfont}{`6}
\DeclareMathSymbol{7}{\mathalpha}{mtoperatorfont}{`7}
\DeclareMathSymbol{8}{\mathalpha}{mtoperatorfont}{`8}
\DeclareMathSymbol{9}{\mathalpha}{mtoperatorfont}{`9}

\DeclareMathSymbol{!}{\mathclose}{mtoperatorfont}{"21}
\DeclareMathSymbol{?}{\mathclose}{mtoperatorfont}{"3F}
\DeclareMathSymbol{*}{\mathalpha}{mtoperatorfont}{"2A}
\DeclareMathSymbol{,}{\mathpunct}{mtoperatorfont}{"2C}
\DeclareMathSymbol{.}{\mathord}{mtoperatorfont}{"2E}
\DeclareMathSymbol{:}{\mathrel}{mtoperatorfont}{"3A} 
\DeclareMathSymbol{;}{\mathpunct}{mtoperatorfont}{"3B}
\edef\mt@minus@sign{\mathord{\expandafter\mathchar\number\mathcode`\-}}
\def\relbar{\mathrel{\smash\mt@minus@sign}}
\def\rightarrowfill{$\m@th\mt@minus@sign\mkern-7mu  %
\cleaders\hbox{$\mkern-2mu\mt@minus@sign\mkern-2mu$}\hfill
  \mkern-7mu\mathord\rightarrow$}
\def\leftarrowfill{$\m@th\mathord\leftarrow\mkern-7mu %
 \cleaders\hbox{$\mkern-2mu\mt@minus@sign\mkern-2mu$}\hfill
  \mkern-7mu\smash\mt@minus@sign$}
\DeclareMathSymbol{-}{\mathbin}{mtoperatorfont}{21}
\DeclareMathSymbol{+}{\mathbin}{mtoperatorfont}{"2B}
\edef\mt@equal@sign{{\expandafter\mathchar\number\mathcode`\=}}
\DeclareRobustCommand\Relbar{\mathrel{\mt@equal@sign}}
\DeclareMathSymbol{=}{\mathrel}{mtoperatorfont}{"3D}
\DeclareMathDelimiter{(}{\mathopen} {mtoperatorfont}{"28}{largesymbols}{"00}
\DeclareMathDelimiter{)}{\mathclose}{mtoperatorfont}{"29}{largesymbols}{"01}
\DeclareMathDelimiter{[}{\mathopen} {mtoperatorfont}{"5B}{largesymbols}{"02}
\DeclareMathDelimiter{]}{\mathclose}{mtoperatorfont}{"5D}{largesymbols}{"03}
\DeclareMathDelimiter{/}{\mathord}{mtoperatorfont}{"2F}{largesymbols}{"0E}
\DeclareMathSymbol{/}{\mathord}{mtoperatorfont}{"2F}
\DeclareMathDelimiter{<}{\mathopen}{mtoperatorfont}{"3C}{largesymbols}{"0A}
\DeclareMathDelimiter{>}{\mathclose}{mtoperatorfont}{"3E}{largesymbols}{"0B}
\DeclareMathSymbol{<}{\mathrel}{mtoperatorfont}{"3C}
\DeclareMathSymbol{>}{\mathrel}{mtoperatorfont}{"3E}
\expandafter\DeclareMathDelimiter\@backslashchar
                        {\mathord}{mtoperatorfont}{"5C}{largesymbols}{"0F}
\DeclareMathDelimiter{\backslash}
    {\mathord}{mtoperatorfont}{"5C}{largesymbols}{"0F}
\DeclareMathSymbol\setminus\mathbin{mtoperatorfont}{"5C}
\DeclareMathSymbol{|}\mathord{mtoperatorfont}{"7C}
\DeclareMathDelimiter{|}{mtoperatorfont}{"7C}{largesymbols}{"0C}
\DeclareMathDelimiter\vert
                 \mathord{mtoperatorfont}{"7C}{largesymbols}{"0C}
\DeclareMathSymbol\mid\mathrel{mtoperatorfont}{"7C}
\DeclareMathDelimiter{\lbrace}
   {\mathopen}{mtoperatorfont}{"7B}{largesymbols}{"08}
\DeclareMathDelimiter{\rbrace}
   {\mathclose}{mtoperatorfont}{"7D}{largesymbols}{"09}

\DeclareMathSizes{9}{9}{7}{5}
\DeclareMathSizes{\@xpt}{\@xpt}{8}{6}
\DeclareMathSizes{\@xipt}{\@xipt}{9}{7}
\DeclareMathSizes{\@xiipt}{\@xiipt}{10}{8}
\DeclareMathSizes{\@xivpt}{\@xivpt}{\@xiipt}{10}
\DeclareMathSizes{\@xviipt}{\@xviipt}{\@xivpt}{\@xiipt}
\DeclareMathSizes{\@xxpt}{\@xxpt}{\@xviipt}{\@xivpt}
\DeclareMathSizes{\@xxvpt}{\@xxvpt}{\@xxpt}{\@xviipt}
\makeatother

\let\oldexists\exists
\renewcommand\exists{\oldexists\,}
\let\oldforall\forall
\renewcommand\forall{\oldforall\,}

\makeatletter
\@for\x:=N,Q,Z,R,C,H\do{%
\edef\act{\noexpand\newcommand{\csname \x\x\endcsname}{\noexpand\mathbf{\x}}}
\act}
\@for\x:=F,O,Z\do{%
\edef\act{\noexpand\newcommand{\csname c\x\endcsname}{\noexpand\mathcal{\x}}}
\act}
\@for\x:=Res\do{%
\edef\act{\noexpand\DeclareMathOperator*{\csname \x\endcsname}{\x}}
\act}
\makeatother

\let\Oh\cO

\renewcommand{\Re}{\mathrm{Re}}
\renewcommand{\Im}{\mathrm{Im}}
\renewcommand\le\leqslant
\renewcommand\ge\geqslant

\let\wh\widehat
\let\bar\overline

\begin{document}

\pdfbookmark[1]{Title Page}{pagedetitre}

\title{On the system of the functions $\frac{\zeta(s)}{(s-\rho)^k}$}
\author{\textsc{Jean-Fran\c{c}ois Burnol}}
\date{\relax}
\maketitle

\vspace{-10\baselineskip}
 \hfill {\small\ladate}\hspace{\the\leftmargin}
\vspace{10\baselineskip}

\begin{abstract}\pdfbookmark[2]{Abstract}{resumei}
  The system of the functions $\smash{\frac{\zeta(s)}{(s-\rho)^k}}$ is
  complete and minimal in a certain sub-Hilbert space of the $L^2$ space of the
  critical line. We study whether it is also hereditarily complete.
\end{abstract}

\vfil

\begingroup\small

\begin{quotation}\pdfbookmark[2]{Thanks}{merci}
  The author thanks the CRM of Barcelona, where this
  work was initiated, for its hospitality, and A.~Baranov,
  Y.~Belov, and the referee for providing useful background information.
\end{quotation}

\vfil

\setlength{\baselineskip}{.75\baselineskip}
\leftskip .4\textwidth \parindent 0pt
\obeylines  
Universit\'e Lille 1   
UFR de Math\'ematiques 
Cit\'e Scientifique M2 
F-59655 Villeneuve d'Ascq 
France
\endgroup

\vfil

\newcommand*{\monheader}{{\normalfont
\ifodd\value{page}%
   {\emph{\dots $\zeta(s)/(s-\rho)^k$}}
  \hfil\thepage
\else 
  \thepage\hfil 
   \emph{On the system of the functions\dots}\fi}}

\markright{\protect\monheader}

\makeatletter
\def\@oddhead{\rightmark}
\makeatother

\vfil
\begin{small}
  to be published in \emph{Complex Analysis and Operator Theory.}
\end{small}
\clearpage

\section{Introduction and statement of the main results}

The functions $\frac{\zeta(s)}{(s-\rho)^k}$, where $\rho$ is a non-trivial
zero of the Riemann zeta function, and $k$ an integer between $1$ and the
multiplicity $m_\rho$ of $\rho$, are square-integrable on the critical line. In
\cite{jftwosystems} I proved that they are a complete and minimal system in
a certain $L\subset L^2(\frac12+i\RR;\frac{|ds|}{2\pi})$.

The Hilbert space $L$ can be characterized as follows: a function $g(s)$
belongs to $L$ if it is the Mellin transform $g(s) = \wh f(s) = \int_0^\infty
f(x) x^{-s}\,dx$ of a square integrable function $f(x)$ on $(0,+\infty)$,
which is constant on $(0,1)$ and such that its cosine transform $\int_0^\infty
2\cos(2\pi xy)f(y)\,dy$ also is constant on $(0,1)$.  I also proved in
\cite{jftwosystems} that the dual (i.e. biorthogonal) system is complete (and minimal, of course)
in $L$.

It is a fact that the Mellin transform $g(s) = \wh f(s)$ of an $f$ satisfying
these support conditions is a meromorphic function in the entire complex
plane, having trivial zeros at $-2$, $-4$, \dots, and at most a pole at
$1$. The entire functions $s(s-1) \pi^{-\frac s2}\Gamma(\frac s2) \wh f(s)$
and $s(s-1) \pi^{-\frac s2}\Gamma(\frac s2) \wh{\cF f}(s)$, where $\cF$ is the
Fourier cosine transform on $L^2(0,+\infty;dx)$, are exchanged by
$s\leftrightarrow1-s$. Evaluating these entire functions or their derivatives
at any given $s$ defines bounded linear forms. I refer to \cite{jftwosystems}
for these and other facts.

A conference talk by {Yurii~Belov} on his joint work with
{Anton~Baranov} \cite{babeimrn} introduced me to the notion of
``hereditary completeness''. 
Under the name of ``strong completeness'', it was defined by
A.~S.~Markus about forty years ago \cite[\textsection3]{markus}:
let $(x_n)_{n\in I}$ be a complete and minimal system in some separable
Hilbert space, and $(y_n)_{n\in I}$ its biorthogonal system (we shall also say
``dual'' for ``biorthogonal''). Let $J\subset I$ and define the vectors $z_n$
by $z_n = x_n$ for $n\in J$ and $z_n = y_n$ for $n\notin J$. If, for all
$J\subset I$, $\cZ_J = (z_n)_{n\in I}$ is a complete system, then $(x_n)$ is
said to be \emph{hereditarily complete}. Equivalently
(\cite[\textsection3]{markus}) the system $(x_n)$ is hereditarily complete if
any vector $x$ is in the closed linear span of the vectors $\ip{y_n,x} x_n$
(scalar products are linear in the second factor).

It follows from the first of these equivalent definitions that a complete and
minimal system which is hereditarily complete has a complete
biorthogonal. Markus constructed in \cite{markus} an example showing that this
necessary condition is not sufficient. I.~N.~Dovbysh and N.~K.~Nikolskii
proposed two simpler, and general, methods leading to such
systems which are not hereditarily complete \cite{dovnik}.

Let us use as index set $I$ the set of all couples $(\rho,k)$ with
$\zeta(\rho) = 0$ (non-trivial zero) and $1\le k \le m_\rho$ and define
$x_{\rho,k} = \frac{\zeta(s)}{(s-\rho)^k}$.  From \cite{jftwosystems} these
vectors are a complete and minimal system and the dual system is also
complete. We study here whether this family is also hereditarily
complete. The technique used to this aim will add some improvements to the
methods from our previous publication \cite{jftwosystems}, but our results are
not complete: hopefully this will stimulate further researches.

To state the result, we need some notations. We do not consider arbitrary
subsets $\Sigma$ of the index set $I$, but only those, which we call
admissible, which are defined in the following manner: $\Sigma = \{ (\rho,k),
1\le k \le k_\Sigma(\rho)\}$ where the function $k_\Sigma: \rho \mapsto
k_\Sigma(\rho)\in\{0,1,\dots,m_\rho\}$ is otherwise arbitrary.  The matrix of
size $m_\rho\times m_\rho$ expressing the dual vectors $y_{\rho,m_\rho}$,
$y_{\rho,m_\rho-1}$, \dots, $y_{\rho,1}$ in terms of the evaluators $g\mapsto
g(\rho)$, $g\mapsto g'(\rho)$, \dots, $g\mapsto g^{(m_{\rho} - 1)}(\rho)$ is
upper-triangular and invertible.\footnote{An explicit formula shall be given
  later.} In particular, requiring that $g$ should be perpendicular to
$y_{\rho,m_\rho}$, $y_{\rho,m_\rho-1}$, \dots, $y_{\rho,k+1}$, is equivalent
to asking that $g(\rho) = g'(\rho) = \dots = g^{(m_\rho-k-1)}(\rho) = 0$
(which is satisfied in particular by the $k$ vectors $g = x_{\rho,1},
x_{\rho,2}, \dots, x_{\rho,k}$).

\begin{theo}\label{th:1}\pdfbookmark[2]{Theorem 1}{theoi}
  Let $\Sigma$ be an admissible subset of the index set $I= \{(\rho,k), 1\le
  k\le m_\rho\}$.  Let $\cZ_\Sigma$ be the system of the vectors $x_{\rho,k}$,
  $1\le k\le k_\Sigma(\rho)$, combined with the $y_{\rho,k}$,
  $k_\Sigma(\rho)<k\le m_\rho$ (or equivalently with the evaluators $g\mapsto
  g^{(j)}(\rho)$, $0\le j < m_\rho - k_\Sigma(\rho)$). The closed linear span
  of $\cZ_\Sigma$ has at most codimension $1$ in $L$.
\end{theo}

{A. Baranov} and {Y. Belov} have studied in a general manner in
\cite{babeimrn} systems of reproducing kernels in some Hilbert spaces
consisting of analytic functions, identifying classes of spaces where the
biorthogonal (we also say ``dual'') system is always
complete, and giving examples where it is not complete. They examine the
question of the hereditary completeness in further work \cite{babe}, and in
particular inside the Paley-Wiener spaces. They have a general ``codimension
at most $1$'' Theorem in this context (and will perhaps in fact exclude, under
general circumstances, the codimension $1$ case). In the present paper, we
don't know whether codimension $1$ is a true possibility or only an indication
of the weaknesses of the techniques we have employed.

To prove \jfautoref{th:1} let's assume to the contrary that there are two
functions perpendicular to the vectors of the system $\cZ_\Sigma$, then there
is one, say $G$, non trivial, and with the additional condition $G(0) = 0$. We
know (see \cite{jftwosystems}) that $g(s)= G(s)/s$ also belongs to $L$. At
each $\rho$ the function $G$, hence also $g$, has order of vanishing at least
equal to $m_\rho - k_\Sigma(\rho)$. And $G$ is perpendicular
to the $\zeta(s)/(s-\rho)^k$, $1\le k \le k_\Sigma(\rho)$.

According to \jfautoref{th:2} (which is stated next) the function $g$
belongs to the closed linear span of the ${\zeta(s)}/{(s-\rho)^k}$, $1\le
k \le m_\rho - m_\rho(g)$, where $m_\rho(g)$ is the multiplicity of $\rho$ as
a (possible) zero of $g$.  But $m_\rho(g)\ge m_\rho -
k_\Sigma(\rho)$, hence $m_\rho - m_\rho(g) \le k_\Sigma(\rho)$. So $G$ is
perpendicular to $g$:
\[ \int_{-\infty}^\infty \frac{|G(\frac12+it)|^2}{\frac12+i t}\,dt = 0\]
Taking the real part we obtain that $G$ vanishes identically, contradiction.

Hence it suffices to prove the following:
\begin{theo}\label{th:2}\pdfbookmark[2]{Theorem 2}{theoii}
  Let $g\in L$ belong to the domain of multiplication by $s$. Then $g$ is in
  the closed linear span of the vectors $\frac{\zeta(s)}{(s-\rho)^k}$, $1\le
  k \le m_\rho - m_\rho(g)$, where $m_\rho(g)$ is the multiplicity of $\rho$ as a
  zero of $g$ (so $0$ if $g(\rho)\neq0$).
\end{theo}

We could obtain the conclusion of \jfautoref{th:2} under weaker hypotheses on
$g$ (in particular under hypotheses which do not exclude from their scope the
functions $\zeta(s)/(s-\rho)$ themselves). But this would add some technical
complications, which anyhow still require some conditions to be imposed
upon the function $g$. The stated formulation thus suffices to our goal here.

\section{Proof of \protect\jfautoref{th:2}}

Let $\phi(x)$ be a smooth function on $(0,\infty)$ with its compact support in
$[\frac1e, e]$, and such that $\wh\phi(\frac12) = \int_0^\infty \phi(x)
x^{-\frac12}\,dx = 1$. We will also use $\psi(x) = \frac1x \phi(\frac1x)$,
which verifies $\wh \psi(s) = \wh\phi(1-s)$. The Mellin transform $\theta(s) =
\wh\phi(s)$ is an entire function which decreases faster than any inverse
power of $|s|$ when $|s|\to\infty$ in any fixed vertical strip of finite width
(follows immediately from integration by parts). 

Note that $\wh\phi(\frac12+it) = \int_{-1}^1 \omega(u)e^{-iut}\,du$, where
$\omega(u) =  \phi(e^u)e^{u/2}$, and that
${\wh\phi(\frac12+i \epsilon t)} = \int_{-\epsilon}^\epsilon
\frac1\epsilon\omega(\frac u\epsilon)e^{-iut}\,du$.  Let $\theta_\epsilon(s) =
\theta(\epsilon(s-\frac12) + \frac12)$. On any compact this converges
uniformly to the constant function $1$ as $\epsilon\to0$, and
$\theta_\epsilon(s)$ is uniformly bounded in $s$ and $0<\epsilon<1$ when $s$
is restricted to a vertical strip of finite width.

Let $g\in L$, $g(s) = \int_0^\infty f(x)x^{-s}\,dx$.  The function
$g_\epsilon(s) = \theta_\epsilon(s) g(s)$ is the Mellin transform of the
multiplicative convolution:
\[ f_\epsilon(x) = \int_{\exp(-\epsilon)}^{\exp(+\epsilon)}
\phi_\epsilon(t) f(\frac xt)\frac{dt}t\] with $\phi_\epsilon(t)\sqrt{t} =
\frac1\epsilon \left( \phi(t^{1/\epsilon})\sqrt{t^{1/\epsilon}}\right)$.

A Mellin transform such as $\int_0^\infty f(x)x^{-s}\,dx$ can also be written
for $s = \frac12+it$ as $\smash{\int_0^\infty
  f(x)x^{1/2}x^{-it}\,\frac{dx}x}$, thus exhibiting it as the additive Fourier
transform of $u\mapsto f(e^u)e^{u/2}$. From this point of view we thus know
that multiplying two Mellin transforms is like the additive convolution of
two functions $\alpha(e^u)e^{u/2}$ and $\beta(e^u)e^{u/2}$ whose result should
be written as a function $\gamma(e^u)e^{u/2}$, thus given by
\[ \gamma(e^u)e^{u/2} = \int_\RR \alpha(e^{u-v})e^{(u-v)/2}\beta(e^v)e^{v/2}
\,dv = \int_\RR \alpha(e^{u-v})\beta(e^v)e^{u/2}\,dv\]\[\implies \gamma(x) =
\int_0^\infty \beta(t)\alpha(\frac xt)\frac{dt}t\] This explains the formula
for the multiplicative convolution $f_\epsilon$. %% =f\ast \phi_\epsilon$.

The function $f_\epsilon$ is constant for $0< x < \exp(-\epsilon)$, and its
Fourier cosine transform also: indeed $\cF(f_\epsilon)$ is the multiplicative
convolution of $\cF(f)$ with $\psi_\epsilon(x) = \frac1x
\phi_\epsilon(\frac1x)$ (the ``Intertwining formula'' of \cite{jfcopoisson}).

I will need a formula for $\theta_\epsilon(s) g(s)$ as a Mellin transform but
for large $\Re(s)$. The expression $\int_0^\infty f_\epsilon(x) x^{-s}\,dx$
needs modification to give an integral which makes sense for $\Re(s)\ge1$,
because of the behavior for $x\to0$. 

Let us first look at pointwise values of
$f_\epsilon(x)$:
\[ |f_\epsilon(x)|^2 x \le \int_{\exp(-\epsilon)}^{\exp(+\epsilon)}
|\phi_\epsilon(x)|^2 dx \int_0^\infty |f(x)|^2\,dx = c \epsilon^{-1}\] for
some constant $c$, so in particular for $\epsilon$ fixed, we have
$f_\epsilon(x) = \Oh(x^{-1/2})$ as $x\to+\infty$. This shows that for any
$\eta>0$, $\smash[b]{\int_\eta^\infty f_\epsilon(x) x^{-s}\,dx}$ makes sense
directly as an analytic function for $\Re(s) > \frac12$. Although we don't
really need it, let us observe that a much better bound can be obtained for
$f_\epsilon(x)$ as $x\to+\infty$. Indeed, with $\cF$ the Fourier cosine
transform, and $\psi_\epsilon(t) = \frac1t \phi_\epsilon(\frac1t)$:
\begin{align*}
  f_\epsilon(x) &= \int_0^\infty \frac1t \phi_\epsilon(\frac xt) f(t)\,dt =
  \int_0^\infty \frac1x \psi_\epsilon(\frac tx) f(t)\,dt\\ 
  &= \int_0^\infty
  \cF(\psi_\epsilon)(xy)\cF(f)(y)\,dy\\
  &= \int_0^1 \cF(\psi_\epsilon)(xy) \beta\,dy + \int_1^\infty
  \cF(\psi_\epsilon)(xy)\cF(f)(y)\,dy\\
  &= -\beta\int_1^\infty \cF(\psi_\epsilon)(xy)\,dy + \int_1^\infty
  \cF(\psi_\epsilon)(xy)\cF(f)(y)\,dy
\end{align*}
Here, $\beta$ is the constant value of $\cF(f)$ on $(0,1)$. Now,
$\cF(\psi_\epsilon)$ is an even function in the Schwartz class, and it follows
then by elementary arguments that $f_\epsilon(x)$ also has Schwartz decrease
as $x\to+\infty$. This is a general phenomenon related to the support
property \cite[\textsection4]{jfcopoisson}. In this manner, we see  that in
fact $\int_\eta^\infty f_\epsilon(x) x^{-s}\,dx$ directly defines an entire
function of $s$, for any $\eta>0$. 

And for $\eta\le\exp(-\epsilon)$,
$f_\epsilon(x)$ is a constant $C(\epsilon)$ on $(0,\eta)$ and we can compute
$\int_0^\eta C(\epsilon) x^{-s}\,dx$ for $\Re(s) <1$, do the analytic
continuation and reexpress it as $- \int_\eta^\infty C(\epsilon) x^{-s}\,dx$
for $\Re(s)>1$. In the end we obtain that a valid representation of
$\theta_\epsilon(s) g(s)$ as an absolutely convergent integral, for $\Re(s) >
1$, is
\[ \int_\eta^\infty (f_\epsilon(x) - C(\epsilon)) x^{-s}\,dx\] where
$\eta$ is chosen $\le\exp(-\epsilon)$. The quantity
$C(\epsilon)$ is also the opposite of the residue of $\theta_\epsilon(s) g(s)$
at $s = 1$, so it is $\theta(\epsilon \frac12 + \frac12)$ times the constant
value $C(0)$ of $f$ on $(0,1)$. We have $\lim_{\epsilon\to 0} C(\epsilon)
= C(0)$, and at any rate this is a bounded quantity. These remarks will serve
later.

The functions $f_\epsilon$ converge to the original $f$ in the $L^2$ sense as
$\epsilon \to 0^+$, but the problem is that the $g_\epsilon$ do not
necessarily belong to $L$: $f_\epsilon$ and $\cF(f_\epsilon)$ are a priori
constant only on $(0,\exp(-\epsilon))$. In the similar computations from my
paper \cite{jftwosystems} this problem was avoided by first replacing $f$ with
a function with stronger support properties, but here we can't do that, at
least we do not see an obvious way to regularize the function $g$ (making it
decrease in the vertical direction) while at the same time maintaining its
vanishing on a certain set of zeros.

There is an a priori (polynomial in vertical strips) upper bound on the growth
of $g(s)$ \cite[Th. 4.8]{jftwosystems}, so $g_\epsilon(s) = \theta_\epsilon(s)
g(s)$ indeed decreases faster than any inverse polynomial when we go to
$\infty$ in any fixed vertical strip of finite width. This allows computing
some contour integrals, with the help of the following theorem:

%% which I believe originates with Valiron, and certainly exists also under
%% stronger forms, but the following will be sufficient for our immediate goals.

    \begin{prop*}[from {\cite[IX.7.]{tit}}]
      There is a real number $A$ and a strictly increasing sequence $T_n>n$
      such that $|\zeta(s)|^{-1} < |s|^A$ on $|\Im(s)|=T_n$, $-1\le \Re(s)
      \le +2$. \end{prop*}

\begin{note}\label{note}
  (taken verbatim from \cite{jftwosystems}) from now on an infinite sum
  $\sum_\rho a(\rho)$ (with complex numbers or functions or Hilbert space
  vectors $a(\rho)$'s indexed by the non-trivial zeros of the Riemann zeta
  function) means $$\lim_{n\to\infty} \sum_{|\Im(\rho)|< T_n} a(\rho)\;,$$ where
  the limit might be, if we are dealing with functions, a pointwise almost
  everywhere limit, or a Hilbert space limit. When we say that the partial
  sums are bounded (as complex numbers, or as Hilbert space vectors) we only
  refer to the partial sums as written above. When we say that the series is
  absolutely convergent it means that we group together the contributions of
  the $\rho$'s with $T_n < |\Im(\rho)| < T_{n+1}$ before evaluating the
  absolute value or Hilbert norm. When building series of residues we write
  sometimes things as if the zeros were all simple: this is just to make the
  notation easier, but no hypothesis is made in this paper on the
  multiplicities $m_\rho$, and the formula used for writing $a(\rho)$ is a
  symbolic representation, valid for a simple zero, of the more complicated
  expression which would apply in case of multiplicity.
\end{note}

Let us follow the method of
\cite[Thm. 5.2]{jftwosystems}, which is to consider a contour integral with
\[ F(s) = \frac{g_\epsilon(s)}{\zeta(s)}\frac{\zeta(Z)}{Z-s}\] where $Z$ is a
fixed parameter. We will mainly be interested by the $Z$'s on the critical
line, but let us take it arbitrarily at this stage (distinct from $1$ and from
the zeros of the Riemann zeta function). We integrate $F(s)$ on the rectangle
with boundary lines 
$|\Re(s) - \frac12| = d$, $|\Im(s)| = T_n$, where $d>\frac12$ is large
enough so that $|\Re(Z) - \frac12|<d$. Letting $n\to\infty$ we obtain:
\[\sum_{\rho}
\frac{g_\epsilon(\rho)}{\zeta^\prime(\rho)}\frac{\zeta(Z)}{Z - \rho}
- g_\epsilon(Z) = \frac{\zeta(Z)}{2\pi}\left(\int_{\Re(s) = \frac12 + d} -
\int_{\Re(s) = \frac12 -
d}\frac{g_\epsilon(s)}{(Z-s)\zeta(s)}|ds|\right)\]

Let us pause to comment on the meaning of
$\frac{g_\epsilon(\rho)}{\zeta^\prime(\rho)}\frac{\zeta(Z)}{Z - \rho}$: as
explained in the \jfautonote{note}, it is a symbolic notation for 
\begin{gather*}
  \Res_{s=\rho} \frac{g_\epsilon(s)}{\zeta(s)}\frac{\zeta(Z)}{Z-s} =
  \sum_{0\le j < m_\rho} \sum_{0\le i \le j} c_{j-i}(\rho)
  \frac{g_\epsilon^{(i)}(\rho)}{i!}\frac{\zeta(Z)}{(Z- \rho)^{m_\rho-j}}\\
\shortintertext{where}
  \frac{(s-\rho)^{m_\rho}}{\zeta(s)} = c_0(\rho) + c_1(\rho) (s-\rho) +
  c_2(\rho) (s-\rho)^2 + \dotsb
\end{gather*}
The linear combination $G\mapsto \sum_{0\le i \le j} c_{j-i}(\rho)
  \frac{G^{(i)}(\rho)}{i!}$ of evaluators, applied to $G(s) = x_{\rho',j'}(s) = 
  \zeta(s)/(s-\rho')^{j'}$ gives $1$ if $\rho'=\rho$ and $j' = m_\rho-j$ and $0$
  otherwise, as can be seen from direct calculation  of $\Res_{s=\rho}
  \frac1{(s-\rho')^{j'} (Z-s)}$; it thus represents the vector
  $y_{\rho,m_\rho-j}$ of the dual system. 

The change of variable $s\mapsto 1-s$ transforms the integral on the line
$\Re(s) = \frac12- d$ into a similar one (where $Z$ is replaced by $1-Z$) on
the line $\Re(s) = \frac12+d$:
\[ \frac{g_\epsilon(1-s)}{\zeta(1-s)} =
\frac{\wh{\cF(f_\epsilon)}(s)}{\zeta(s)}\] As we have already mentioned that
$\cF(f_\epsilon)$ is the multiplicative convolution of $\cF(f)$ by
$\psi_\epsilon$, all our future arguments and bounds for the integral
initially already defined on the line $\Re(s) = \frac12+d$ would apply
similarly to the integral initially on $\Re(s) = \frac12-d$.

On the line $\Re(s) = \frac12+ d$,  $\zeta(s)^{-1}$ can be replaced with the
absolutely convergent expression $\sum_{k\ge 1} \mu(k) k^{-s}$, which
allows termwise integration.
Let us check that for $\epsilon\le\log2$ all the contributions with $k\ge2$
vanish.  For this we write $g_\epsilon(s) k^{-s} = \int_0^\infty \frac1k
f_\epsilon(\frac xk) x^{-s} \,dx$. From previous discussion we know that
the correct formula when $\Re(s)>1$ is:
\[ \int_\eta^\infty \left( \frac1k f_\epsilon(\frac xk) - \frac1k
  C(\epsilon)\right) x^{-s} \,dx\] with some $\eta \le k \exp(-\epsilon)$. For
$k\ge2$ and $\epsilon \le \log 2$ we can take $\eta = 1$ in this formula. We
want to evaluate
\[ \frac{1}{2\pi} \int_{\Re(s) = \frac12 + d} \frac{g_\epsilon(s)k^{-s}}{Z-s}|ds|\]
as an application of Plancherel theorem.\footnote{It is also possible to shift
  the contour of integration to the right to show that it vanishes for $k\ge2$
  and $\epsilon\le\log2$.} So we compute the c.c. (complex
conjugate):
\[ c.c.(Z - s) = \bar Z - (1+2d - s) = - (1+2d - \bar Z - s)\]
With $w = 1+2d - \bar Z$, there holds $\Re(w) > \frac12 + d = \Re(s)$, so 
\[ (w - s)^{-1} = \int_0^1 x^{w-1 - s} \,dx = \int_0^1 x^{d - \bar Z} x^{-\frac12
  - i \Im(s)}\,dx \]
On the other hand:
\[ g_\epsilon(s) k^{-s} = \int_\eta^\infty \left( \frac1k
  f_\epsilon(\frac xk) - \frac1k C(\epsilon)\right) x^{-d} x^{-\frac12- i
  \Im(s)} \,dx\] So, by the Plancherel formula:
\[ \frac{1}{2\pi} \int_{\Re(s) = \frac12 + d}
\frac{g_\epsilon(s)k^{-s}}{Z-s}|ds| = - \int_{\min(\eta,1)}^1 x^{-Z} \left( \frac1k
  f_\epsilon(\frac xk) - \frac1k C(\epsilon)\right) dx\] For $k\ge2$ (and
$\epsilon \le \log 2$) we can take $\eta=1$ and this vanishes.

So we have the representation, for each given fixed $Z$ (not $1$ and not a
zero of the Riemann zeta function):
\[\sum_{\rho} \frac{g_\epsilon(\rho)}{\zeta^\prime(\rho)}\frac{\zeta(Z)}{Z -
    \rho} - g_\epsilon(Z) = {\underbrace{\frac{\zeta(Z)}{2\pi}\int_{\Re(s) =
      \frac12 + d} \frac{g_\epsilon(s)}{Z-s}|ds|}_{A_\epsilon(Z)}} +
  B_\epsilon(Z) \]
\[
  A_\epsilon(Z) = - \zeta(Z)\int_\eta^1 x^{-Z} \left( f_\epsilon(x) -
    C(\epsilon)\right) dx\qquad \eta = \exp(-\epsilon)
\]

The convergence of the series taken over the zeros of the Riemann zeta
function (and with the meaning from the \jfautonote{note}) has so far only
been proven pointwise. The second half of \cite[Proof of 5.2]{jftwosystems}
gives, on page 80, arguments to establish that the series of functions of $Z$
indexed by the zeros of the Riemann zeta function (and their multiplicities)
is an absolutely convergent one in the sense of the $L^2$-norm (and with the
meaning from the \jfautonote{note} above). We do not repeat the arguments
which can be applied here identically. As a corollary the sum $A_\epsilon(Z) +
B_\epsilon(Z)$ is square-integrable on the critical line, a fact which is seen
directly from $\int_\eta^1 x^{-Z} \left( f_\epsilon(x) - C(\epsilon)\right)dx
= \Oh(\frac1{1+|Z|})$ for $\Re(Z)$ bounded, obtained by an integration by
parts, as $f_\epsilon$ is smooth.  But we would also like to examine, as this
would complete the proof of Theorem \ref{th:2}, if the $L^2$-norm of
$A_\epsilon(Z) + B_\epsilon(Z)$ goes to zero as $\epsilon\to0$; this is where
we will use the hypothesis that $sg(s)$ also belongs to $L$.

As an aside, for a fixed $Z$ we can show without hypothesis that
$A_\epsilon(Z)$ goes to zero. We already mentioned that $C(\epsilon)$ was
bounded, and we estimated pointwise $|f_\epsilon(x)| \le c
(x\epsilon)^{-\frac12}$ for some constant $c$. As we integrate over
$x$ in the interval from $e^{-\epsilon}$ to $1$, this gives $A_\epsilon(Z) =
\zeta(Z) \Oh(\epsilon^{\frac12})$, uniformly in $Z$ for $\Re(Z)$ bounded.

We now bound $A_\epsilon(Z)$ otherwise. As we are mainly interested in $\Re(Z)
= \frac12$, we will from now on take $d= 1$. By the Plancherel
argument, or by a shift of the line of integration towards $+\infty$:
\[ \frac{1}{2\pi}\int_{\Re(s) = \frac32} 
\frac{g(s)}{Z-s}|ds| = 0\;.\]
\[ A_\epsilon(Z) = \frac{\zeta(Z)}{2\pi}\int_{\Re(s) = \frac32} 
\frac{g_\epsilon(s) - g(s)}{Z-s}|ds|\]
Writing $\frac1{Z-s} = \frac1Z + \frac{s}{Z(Z-s)}$ and using Cauchy-Schwarz:
\[ |A_\epsilon(Z)|\le \frac{|\zeta(Z)|}{|Z|} 
\int_{\Re(s) = \frac32}  {|s||g_\epsilon(s) - g(s)|}\left(\frac1{|s|} + \frac1{|Z-s|}\right)\frac{|ds|}{2\pi}\]
\[ |A_\epsilon(Z)|\le \frac{|\zeta(Z)|}{|Z|} \sqrt{\int_{\Re(s) = \frac32}
  {|s|^2|g_\epsilon(s) - g(s)|^2}\frac{|ds|}{2\pi}}\left(\sqrt{\frac13} +
  \sqrt{\frac12}\right)\] The last remaining integral does not depend on $Z$
but is a numerical quantity depending on $\epsilon$. It goes to zero as
$\epsilon\to0$ from the Lebesgue dominated convergence theorem. We silently
used that $sg(s)$ was square-integrable on the line $\Re(s) = \frac32$. But
this is clear as, by hypothesis, $sg(s) = \frac{C}{s-1} + k(s)$ with some $k$
in the Hardy-space of the half-plane $\Re(s) > \frac12$.

Combining the results obtained we conclude that $g(Z)$ can be arbitrarily well
approximated in $L^2$-norm by a finite linear combination of the
$\zeta(Z)/(Z-\rho)^k$ where only those $k$ between $1$ and $m_\rho -
m_\rho(g)$ (inclusive)
appear, which is the statement of \jfautoref{th:2}.

\section{The completeness of the evaluators without Kre\u\i n's theorem}

In \cite{jftwosystems} I proved that the evaluators associated with the zeros of
the Riemann zeta function were complete: i.e. if an element $g$ in $L$ is such
that $g(s)/\zeta(s)$ is entire, then $g$ is the zero function. I used a Theorem
of Kre\u\i n on entire functions in the Cartwright class. 

A more elementary proof can now be given. Again with $g_\epsilon(s)$ being defined
as $\theta(\epsilon(s-\frac12)+\frac12)g(s)$, in the evaluation of the contour
integral built with $\frac{g_\epsilon(s)}{\zeta(s)}\frac{\zeta(Z)}{Z-s}$
(where $Z$ is again a parameter distinct from $1$ and from the zeros of the
zeta function) the only singularity is now at $s = Z$, and we obtain the
formula:
\[ 
-g_\epsilon(Z) = \frac{\zeta(Z)}{2\pi}\left(\int_{\Re(s) = \frac12 + d} -
\int_{\Re(s) = \frac12 -
d}\frac{g_\epsilon(s)}{(Z-s)\zeta(s)}|ds|\right) = A_\epsilon(Z) + B_\epsilon(Z)\]
We can as well take $\Re(Z) = \frac12$ and $d = 1$. But we have argued already
that for fixed $Z$ there hold (under no additional hypothesis on $g$) the
pointwise limits $A_\epsilon(Z) \to 0$, $B_\epsilon(Z)\to 0$, for $\epsilon\to
0$. This proves that $g$ is the zero function.

The same argument would show that the only functions in $L$ which vanish (with
at least the same multiplicities) on
all but perhaps finitely many zeros of the Riemann zeta function are the
finite linear combinations of the functions
$\frac{\zeta(s)}{(s-\rho)^k}$. Indeed the sum of the residues being now finite,
there is no problem with taking the limit $\epsilon\to0$ to obtain a
pointwise identity, which suffices for the conclusion.

This gives examples of mixed systems being complete, but I must leave open the
question whether codimension $1$ can really happen for some other kind of
combined system.

\end{document}